\documentclass[11pt]{amsart}
\usepackage{amsmath}
\usepackage{amssymb}
\usepackage{amsthm}
\usepackage[pdftex]{graphicx}
\newcommand{\ignore}[1]{\relax}

\newcommand{\C}{\mathbb C}
\newcommand{\R}{\mathbb R}
\newcommand{\Z}{\mathbb Z}

\newcommand{\Q}{\mathbb Q}

\newcommand{\Arg}{\operatorname{Arg}}

\newtheorem{lem}{Lemma}

\theoremstyle{definition}

\newtheorem{exa}[lem]{Example}

\theoremstyle{remark}

\newtheorem{rmk}[lem]{Remark}

\newcommand{\dd}{\partial}

\newcommand{\cp}{{\mathbb C}{\mathbb P}}
\newcommand{\tp}{{\mathbb T}{\mathbb P}}
\newcommand{\rp}{{\mathbb R}{\mathbb P}}

\newcommand{\Log}{\operatorname{Log}}

\newcommand{\T}{\mathbb{T}}

\renewcommand{\setminus}{\smallsetminus}

\begin{document}

\title{Geometry in the tropical limit}
\author{I. Itenberg, G. Mikhalkin}
\address{Universit\'e Pierre et Marie Curie and
Institut Universitaire de France\\
Institut de Math\'ematiques de Jussieu\\
4 place Jussieu\\
75005 Paris, France}
\email{itenberg@math.jussieu.fr}
\address{Universit\'e de Gen\`eve\\
Math\'ematiques, villa Battelle\\
7, route de Drize\\
1227 Carouge, Switzerland}
\email{grigory.mikhalkin@unige.ch}
\begin{abstract}
Complex algebraic
varieties
become easy piecewise-linear objects after passing
to the so-called tropical limit. Geometry of these limiting objects is known as tropical geometry.
In this short survey we take a look at motivation and intuition behind this limit and consider
a few simple examples of correspondence principle between classical and tropical geometries.
\end{abstract}

\thanks{This paper is an extended version of the contribution to the 2011 Jahrbuch of the Max-Planck-Gesellschaft
http://www.mpg.de/1368537.
The authors thank the Max-Planck-Institute f\"ur Mathematik in Bonn for its hospitality on a number of occasions.
Research is supported in part by
ANR-09-BLAN-0039-01 grant of {\it Agence Nationale de la
Recherche}
(I.I.),
by the TROPGEO project of the European Research Council and
by the Swiss National Science Foundation  grants 125070 and 126817 (G.M.)
as well as FRG 0854989/0854977 grants of the National Science Foundation USA (I.I. and G.M.)
}
\maketitle

\section{Introduction} Algebraic geometry studies geometric objects associated
to polynomial equations. Such equations make sense over any choice of coefficients
as long as we can add and multiply them (subject to the usual commutativity, associativity
and distribution law). Quite often one chooses an algebraically  closed field, such as
the field $\C$ of complex numbers. However one can consider algebraic geometry not
only over other fields, such as the field $\R$ of real numbers, or the field $\Q$ of rational
numbers, but also with coefficients that do not form a field or for that matter not even a ring.

In this survey we look at what happens if we take the so-called tropical numbers $\T$ for coefficients.
Set-theoretically we may take $\T=\R\cup\{-\infty\}$ and enhance it with $\max$ for addition
and $+$ for multiplication. The result is not a field as $\max$ is an idempotent operation
and does not admit an inverse. Nevertheless there are meaningful geometric objects,
called tropical varieties, associated to tropical polynomials.

Tropical arithmetic operations appear as a certain limiting case  of classical  additions and multiplications.
Given two expressions $\alpha t^{a}+o(t^{a})$ and $\beta t^{b} + o(t^b)$, which are monomial in $t$
with $o$-small precision, $t\to +\infty$,
their sum has the form $\gamma t^{\max\{a,b\}} + o(t^{\max\{a,b\}})$ while
their product has the form $\alpha\beta t^{a+b} + o(t^{a+b})$.
If  $\alpha\beta\neq 0$ and $\alpha+\beta\neq 0$ then  rough asymptotic behavior of the four expressions
is determined by $a$, $b$, $\max\{a,b\}$, $a+b$, respectively, resulting in appearance of tropical operations.
In their turn tropical varieties may be presented as results of collapse of complex algebraic varieties and
through this can be viewed as limiting complex objects.

One can meet such type of limit quite often in various areas of science, e.g. Quantum Mechanics and Thermodynamics.
We start this survey by reviewing how this limit appears there, particularly, the relevance of complex numbers
in quantum formalism as well as thermodynamical interpretation of pre-tropical (subtropical) addition.

\section{Complex numbers and their quantum-mechanical motivation}\label{dequantization}
In Mathematics complex numbers are traditionally considered as the most natural choice of coefficients.
For most mathematicians these are the easiest imaginable type of numbers to work with. Unlike the situation with
the real numbers, any polynomial equation with complex coefficients has solutions. Yet complex numbers
are easy to visualize by thinking of them as points on the 2-plane.
\begin{figure}[h]
\includegraphics[height=25mm]{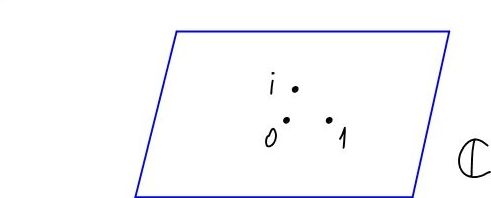}
\caption{\label{cplane} Complex numbers}
\end{figure}

But is such a viewpoint actually supported by non-mathematical considerations? Of course as of today
we have not seen the appearance
of numbers like $-\frac{\sqrt3}{2}+\frac{i}{2}$ in Geography or even in Biology.
Nevertheless, since at least the middle of the XIXth century (ever since the discovery of Electromagnetism)
the complex numbers
are
a quintessential tool in Physics. Namely, a complex number $z=re^{i\phi}$ possesses
the {\em phase} $\phi$.
Alternating current (that is available to us from a household electric socket) can be described
by a complex number whose phase changes in time
(e.g. to make the frequency of 50Hz the phase has to increase
by $100\pi$,
{\it i.e.},
by $50$ full circles around the origin
in the complex plane $\C$ every second).

In the beginning of the XXth century these ideas were greatly advanced in quantum physics. According to
Schr\"odinger each physical particle can be thought of as a probabilistic distribution of its possible coordinate values
plus the {\em choice of phase} at every point of the physical space. The motion of the particle is described not
only by change of its distribution but also by change of its phase with time.
In particular, the celebrated formula $E=\hbar\omega$ of Max Planck expresses the energy of a particle
(in a stationary state) through the frequency of its phase change.

Let us recall how the Planck formula can be interpreted in terms of Schr\"odinger's wave function.
In classical Mechanics we think of a particle in the 3-space $\R^3$ as a point $x(t)\in\R^3$.
This point changes with time $t\in\R^3$.
Once we pass to a (non-relativistic) Quantum mechanics viewpoint we
may think of a particle as a complex-valued function
$$\psi(t):\R^3\to\C$$
subject to the condition $\int\limits_{\R^3} |\psi(t)|^2 =1$.
Thanks to this condition, the real-valued function $|\psi(t)|^2$ is a time-varying probability distribution in $\R^3$.
It can be interpreted as probability to meet our particle at a specified position at time $t$.

The argument (phase) of $\psi$ does not have an immediate physical meaning.
The change of $\Arg(\psi(t))$ in space affects the gradient $\nabla\psi(t)$
(which can be interpreted as the momentum operator once multiplied by $-i{\hbar}$).
The change of $\Arg(\psi(t))$ in time is governed by the Schr\"odinger equation.
$$i{\hbar}\frac{d\psi(t)}{dt}=H\psi(t),$$
where $H$ is the Hamiltonian operator acting on the space of all complex-valued $L^2$-functions in $\R^3$.
Eigenvalues of $H$ are called the energy spectrum, the corresponding
eigenfunctions are stationary states.
If $\psi(0)$ is a stationary state then $\psi(t)$ is also a stationary state corresponding to the
same energy level $E$. Furthermore, from the Schr\"odinger equation we have
$$\psi(t)=e^{-i\frac{E}{\hbar}t}\psi(0).$$
In the right-hand side of this equation the factor $\frac{E}{\hbar}$ corresponds to the frequence $\omega$ of
the phase-change of $\psi$ at every point of $\R^3$, so $\omega=\frac{E}{\hbar}$ as in Planck's formula.
The Planck constant $\hbar$ is thus the universal constant that converts frequency units into energy units,
$$1\ Hertz \sim 7 \times 10^{-34}\ Joules.$$

We have a similar situation with the change of frequency in space.
If the phase frequency
is very high and changes rather slowly,
people speak of quasiclassical motion of a quantum particle.
In such cases we may ignore the phase. E.g. we think of the presence of electricity in the household
socket even though at some (rather frequent) moments the real part of the phase vanishes.
Quasiclassical approximation is used to relate classical and quantum mechanics and
provide intuition for the so-called correspondence principle in quantum mechanics.

\section{Can we forget the phase in a complex number?}
To forget the phase $\phi$ in $z=re^{i\phi}$,
it suffices to consider the absolute value $|z|=r$ instead of $z$.
But our goal is to get rid of $\phi$ while keeping basic features of the complex numbers.
In particular, we would like
to keep our ability to add and multiply the numbers regardless of their phase.
To do this we have to pass to a certain limit, called the {\em tropical limit}, introducing a large positive
parameter $t>>1$ which will tend to $+\infty$.

Consider the base $t$ logarithm map $$\Log_t:\C\to\R\cup\{-\infty\}=\T$$ of the absolute value, $z\mapsto\log_t|z|$.
The target $\R\cup\{-\infty\}$ of this map is usually denoted
by $\T$. The  elements of this set are called
{\em tropical numbers}.
We may use the map $\Log_t$ to induce the addition and multiplication operations on $\T$ from $\C$.
\begin{figure}[h]
\includegraphics[height=35mm]{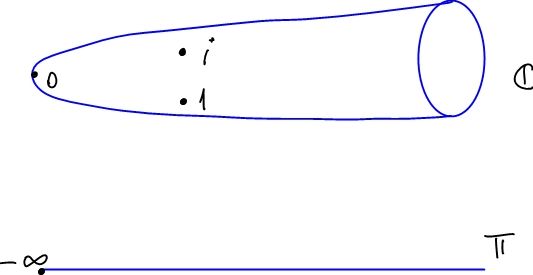}
\caption{\label{tplane} Collapse of complex numbers to tropical numbers}
\end{figure}

It is very easy to define the product of two tropical numbers $x, y\in\T$.
Their inverse images
under $\Log_t$ are $\alpha t^x$ and $\beta t^y$,
$\alpha, \beta\in\C$, $|\alpha|=|\beta|=1$.
We get the induced product of $x$ and $y$
equal to $\log_t|\alpha t^x\beta t^y|=x+y$.
We see that it
depends neither on $\alpha$ and $\beta$
nor on the parameter $t$.
The operation
$``xy"=x+y$  is called {\em the tropical product} of $x,y\in\T$.

The induced sum is
\begin{equation}\label{plus-t}
\log_t|\alpha t^x + \beta t^y|.\
\end{equation}
For a given $t$ it depends on $\alpha$ and $\beta$.
Suppose, say that $x\ge y$, so that $x=\max\{x,y\}$. By the triangle  inequality then we get
$t^x(1-t^{y-x})=t^x-t^y\le|\alpha t^x + \beta t^y|\le t^x+t^y\le 2t^x$.
Taking $\log_t$ of the upper bound we get $x+\log_t 2$ which tends to $x$ when $t\to+\infty$.
Taking $\log_t$ of the lower bound
we get $x+\log_t|1-t^{y-x}|$. This tends to $x=\max\{x,y\}$
if $x>y$, but it is $-\infty$ if $x=y$.
The operation $``x+y"=\max\{x,y\}$ is called {\em the tropical sum} of $x,y\in\T$.
We see that it is the genuine limit of the induced operation from $\C$ whenever $x\neq y$ and it
is the upper limit of such operation if $x=y$.

Similarity between passing to the tropical limit and doing the procedure
inverse to quantization
was noted by Maslov. He and his school have established
a number of theorems in analysis
that correspond to each other under this procedure,
see \cite{Maslov}. Relations with quantum-mechanical notions can be found
in Litvinov's paper \cite{Litvinov}.
This procedure
is also known as {\em Maslov's dequantization}.
To relate it with the quasiclassical limit
one has to set $t=e^{\frac1{\hbar}}$,
so that indeed $t\to+\infty$ is equivalent to $\hbar\to 0$.
Viro observed that his patchworking technique~\cite{V-patchworking}
(the most powerful technique known for construction
of real algebraic varieties) can be obtained through
a quantization inverse to
Maslov's dequantization
of the complex plane, see \cite{V}.

As an aside note we also get a very interesting geometrical situation
if we forget the phase without
passing to the tropical limit. Namely we may consider images
of subvarieties of $(\C^\times)^n$
under the coordinatewise $\Log_t$ map for a finite $t>1$.
This geometry was introduced by
Gelfand, Kapranov and Zelevinsky \cite{GKZ}.
The resulting images in $\R^n$ are called {\em amoebas}.

\section{Tropical addition and zero-temperature limit in thermodynamics}
If we look more closely at the relation between tropical limit and the quasiclassical
limit in quantum mechanics we may notice a twist by $i$. E.g. to get
a rough (leading order in $\hbar$) approximation for the phase of the Schr\"odinger wave function $\psi$ we write
$$\psi(x)=a(x)e^{\frac{iS(x)}{\hbar}},$$
where $S$ is the classical action functional, see \cite{LL3} and $a:\R^3\to\R$ is some real-valued function
 (alternatively we can write $S(x)=\hbar\Arg(\psi(x))$
to express action through the argument of the wave function). Appearance of $i$ in front of
the real-valued function $S$ is notable and is the subject of the famous {\em Wick rotation by $i$}
relating quantum mechanics and thermodynamics (introducing among other things the concept of imaginary
time, much celebrated in popular culture).

This makes thermodynamics another major physical context where
expressions such as \eqref{plus-t} appear naturally (in a sense even more naturally
than in quantum mechanics as rotation by $i$ is no longer needed).
If we set $\alpha=\beta=1$ then \eqref{plus-t} can be viewed
as an addition operation
\begin{equation}\label{oplust}
x\oplus_t y=\log_t (t^x + t^y),
\end{equation}
$x,y\in\R$, parameterized by a positive number $t\neq 1$.
For any such $t$ this operation and the tropical
multiplication $``xy"=x+y$
satisfy
the distribution law
$$\displaylines{
``(x\oplus_t y)z" = \log_t (t^x + t^y)+z = \log_t ((t^x+t^y)t^z) \cr
= \log_t(t^{x+z}+t^{y+z}) = ``xz"\oplus_t ``yz".
}
$$
When $t\to+\infty$ the limit $x\oplus_{\infty}y=\max\{x,y\}$ is the tropical addition.
When $t\to 0$ the limit $x\oplus_0y=\min\{x,y\}$ can be identified with the tropical addition
by the isomorphism $\R\to\R$, $x\mapsto -x$ that preserves tropical multiplication $``xy"$.
Thus $\min\{x,y\}$ can also be viewed as the tropical addition for a different, but isomorphic choice of the model
of tropical arithmetic operations on $\R$.
\footnote{Sometimes in tropical literature $\min$ is chosen as the model for tropical addition on $\R$.}
For the connection to thermodynamics it is more convenient to use this alternative $\min$-model
of tropical addition.

Thus in both limiting cases $t=0$ and $t=+\infty$ we get tropical addition (in $\max$ and $\min$-model).
We call the arithmetic operation \eqref{oplust} for finite positive $t\neq 1$ {\em subtropical $t$-addition}.
Clearly the subtropical addition \eqref{oplust} is an increasing function of $t$.

Starting from the time of steam engine, most of
the machines
that work for us now
are based on one of many possible thermodynamical cycles (e.g. the Otto or Diesel cycles).
There is the working body (in the simplest case we may assume that it is ideal gas in a box) that changes
its state while performing work (outside this system), but at the end of the cycle returns to its initial state.

Let us remind some basic thermodynamical concepts in their simplest, quantum non-relativistic form.
The working body in our thermodynamical system is assumed to be a vessel with ideal quantum Boltzmann gas.
This system has the energy spectrum $E_j$, $j=0,\dots,+\infty$, that is an increasing infinite
sequence, each $E_j$ corresponding to the $j$th stationary state
of the system.

As we assume our gas to be ideal, its particles do not interact with each other (furthermore,
we assume it to be sparse, so that the average number of particles in any given state is
much less than 1, so that we may even neglect the exchange interaction).
Thus the energy of the system is simply the sum of the energies of the individual particles.
Each quantum particle can be in one of
infinitely many stationary
states
(or in a mixed state).

These states are characterized by their energy $\epsilon_j$ and the numbers $E_j$ are obtained
as the sum of possible values of $\epsilon_j$ over the number $N$ of particles and practically almost always
we may assume that all $N$ values for $\epsilon_j$ are different.
The sequence $\epsilon_j$, $j=0,\dots,+\infty$ is determined by such things as the type of gas and the shape of
the ambient vessel
(to find it mathematically we have to solve the corresponding Schr\"odinger equation).

The state of our thermodynamical system is a probabilistic measure on the stationary states of the system
(a countable set in our case). According to the Gibbs law, if we assume our system to be in thermodynamical
equilibrium,
then the probability of the $j$th state is proportional to the weights $e^{-\frac{E_j}{T}}$, where $T>0$ is a parameter
called the {\em temperature} of the system, see \cite{LL5}. \footnote{For simplicity here we measure temperature in the energy units,
otherwise we need to multiply $T_j$ by the
Boltzmann constant converting temperature into energy, $k\sim\  1.4 \times 10^{-23}Joule/Kelvin$.}

The Helmholtz free energy $F$ is $T$ times the logarithm of the partition function associated to these weights:
\begin{displaymath}
F=-T\log(\sum\limits_{j=0}^{\infty} e^{-\frac{E_j}{T}}).
\end{displaymath}
It can be shown that increment of $F$ during an {\em isothermal} process ({\it i.e.}, a process perhaps changing
the energy of the stationary state of the working body, but keeping the temperature constant) equals to the amount
of mechanical work performed on our working body (so that the increment is negative if the working
body performs work).

Note that if we set $t=e^{-\frac1T}$ then
\begin{equation}\label{helmholtz}
F= E_0\oplus_t E_1\oplus_t\dots\oplus_t E_j \oplus_t\dots,
\end{equation}
{\it i.e.},
nothing else but the subtropical $t$-sum of the energies $E_j$ of the stationary states of the system
with parameter $t=e^{-\frac1T}$. Note that the $t\to 0$ limit corresponds to the $T\to 0$ limit,
{\it i.e.}, the tropical limit corresponds to the zero-temperature limit.

Thermodynamical motivations entered considerations in geometry on a number of occasions.
Kenyon, Okounkov and Sheffield \cite{KOS} succeedded in exhibiting amoebas of plane
complex algebraic curves
as limiting objects associated to a certain statistical model (the {\em dimer model}) enhanced with
Gibbs measures. Corresponding geometric objects at zero temperature there can be interpreted
as tropical curves in the plane.

A very inspiring thermodynamical interpretation of toric geometry and, in particular, amoebas was
recently suggested by Kapranov \cite{Ka}.
Recall (see \cite{GKZ}) that an amoeba is
an image of a variety in $(\C^\times)^n$ under the map $\Log: (\C^\times)^n \to \R^n$
defined coordinatewise by the logarithm of the absolute value.
Suppose that $$A\subset\Z^n\subset\R^n$$ is a finite set which we can interpret as the
set of stationary states of a thermodynamical system. Each linear function in $\R^n$
associates energy levels to elements of $A$. In this sense the embedding $A\subset \R^n$
can be thought of as specifying $n$ commuting Hamiltonians (or a vector Hamiltonian).
A point $x$ in the convex hull $\Delta$ of $A$ can be interpreted as a probability measure (convex
linear combination) to be in one of the stationary states. No matter how big is $A$
there is a unique way to present $x=(x_1,\dots,x_n)$ so that the Gibbs law will hold for all $n$ Hamiltonians.
This gives us $n$ temperatures $T_j$, $j=1,\dots,n$, and accordingly
$n$ Boltzmann parameters $\beta_j=\frac1{T_j}$
that serve as coordinates in $\R^n$ viewed as the target space of the map $\Log$.

This interpretation allows to identify canonically the interior of $\Delta$ with $\R^n$
associating the inverse temperatures $\frac1{T_j}$ to the only thermodynamically stable
state with the average energy $x_j$, $j=1,\dots,n$.
Recall that Viro's patchworking \cite{V-patchworking} can be thought of as
gluing real algebraic curves defined on faces of $\Delta$ of $A$ by moving
them slightly off $\dd\Delta$. According to \cite{Ka} this can be interpreted
as passing from zero temperature ($\frac{1}{T}=\infty$) to non-zero low temperature.

\begin{figure}[h]
\includegraphics[height=50mm]{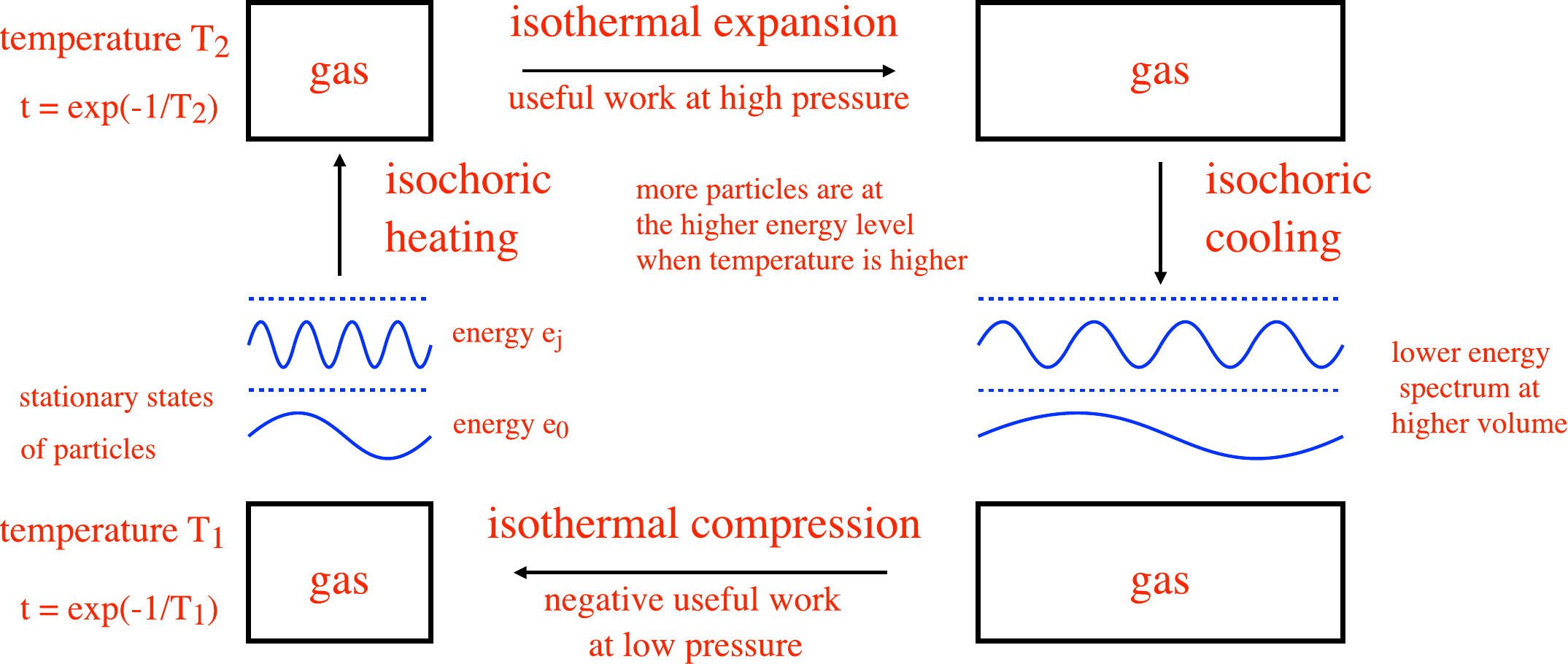}
\caption{\label{stirling} Stirling thermodynamical cycle}
\end{figure}

Here we would like to consider a much simpler example of such a correspondence
based on the so-called {\em Stirling cycle} in thermodynamics.
The Stirling cycle consists of four steps, see Figure \ref{stirling}.
At step I the vessel with gas is heated from a temperature $T_1$ to a temperature $T_2>T_1$
keeping the volume of gas in the vessel fixed (the isochoric heating).
At step II the gas performs work over an exterior system: the gas is allowed to expand isothermally
at the temperature $T_2$ so that it can make useful work, e.g. to move
the pistons in our engine at high pressure.
At step III the gas is isochorically cooled back to the temperature $T_1$.
At step IV the gas is isothermally compressed to its initial state.

Note that some work is performed on the gas at step IV,
{\it i.e.}, in a sense the gas is performing
a negative useful work. However since $T_1<T_2$ the gas pressure will be lower and
the amount of work needed to perform in step IV is less than the useful work performed by our gas in step II.
Thus the useful mechanical work done during the Stirling cycle
is equal to the amount of free energy lost in step II minus the amount of free energy gained in step IV.

In step I the free energy $F$ increases since the subtropical $t$-addition increases as $t$ grows with the temperature $T$,
in step III it decreases. Thus the amount of useful work during the Stirling cycle is bounded from above
by the differences of the subtropical $t$-sums \eqref{helmholtz} at $t=e^{-\frac1{T_2}}$ and $t=e^{-\frac1{T_1}}$.

In the tropical limit $T=0$ (we have $t=e^{-\frac1T}$ as well) the free energy \eqref{helmholtz} just equals to the
energy $E_0$ of the ground state of the gas.

\section{Some tropical varieties and examples of correspondence principle}
The tropical operations described above
give rise to certain meaningful geometric
objects, namely, the {\em tropical varieties}.
From the topological point of view,
the tropical varieties are piecewise-linear polyhedral complexes
equipped with a particular geometric structure
which can be seen as the degeneration
in the tropical limit
of the complex structure
of an algebraic variety.

It is especially easy
to describe tropical varieties in dimension~$1$,
{\it i.e.},  {\em tropical curves}.
Consider, first, tropical curves in
the tropical
affine space $\T^n = (\R \cup \{-\infty\})^n$.
Such a tropical
curve can be obtained as the limit
of the images of some complex algebraic curves
$C_t \subset \C^n$ under the map $\Log_t$, $t \to +\infty$.
The limiting objects
are finite graphs with straight edges (some of them
going to infinity); each edge
of the graph is of rational slope, and
a certain balancing (or ``zero-tension") condition is satisfied
at each vertex of the graph.

There are two natural ways to describe plane curves:
by equation and by parametrization.
Thus, to describe a tropical curve in~$\T^2$, we can either
provide a tropical polynomial defining the curve, or
represent the curve as
the image of an abstract tropical
curve under a tropical map.

A {\em tropical polynomial} in~$\T^2$ (in two variables~$x$ and~$y$)
is
an expression of the following form:
$$``\sum\limits_{(i,j)\in V} a_{ij}x^iy^j"=   \max_{(i,j) \in V}\{a_{ij} + ix + jy\},$$
where~$V \subset \Z^2$ is a finite set of points with
non-negative coordinates and the coefficients $a_{i, j}$
are tropical numbers. The tropical curve
defined by such a polynomial is given by the {\em corner locus}
of the polynomial,
{\em i.e.}, the set of points in $\T^2$, where the function
$$f: (x, y) \mapsto \max_{(i,j) \in V}\{a_{ij} + ix + jy\}$$
is not locally affine-linear. In other words, the corner locus
is the image of ``corners"
of the graph of~$f$ under the vertical projection.

The corner locus of $f$ is composed of intervals and rays in $\R^2$
that form edges of a piecewise-linear graph $\Gamma_f \subset \R^2$.
Each edge $E$ is determined by a choice of two monomials
$``a_{i_1j_1}x^{i_1}y^{j_1}"$ and $``a_{i_2j_2}x^{i_2}y^{j_2}"$ in $f$,
and consists of points where these two monomials coincide and are
larger than other monomials. The GCD of $i_1-i_2$ and $j_1-j_2$
is called the {\em weight} of $E$ and is denoted by $w(E)$.

The edge $E$ is parallel to
the vector $(j_1-j_2,i_2-i_1)=u(E)$ defined by $E$ up to sign (as our two
monomials are not ordered). However, once an endpoint $v$ 
of $E$ is chosen 
(a vertex of graph $\Gamma_f$ adjacent to the edge $E$) then we define
$u(E)$ to be directed away from $v$.
At each vertex $v\in\Gamma_f$ we have the balancing condition:
\begin{equation}\label{balancing}
\sum\limits_{E} u(E)=0,
\end{equation}
where the sum is taken over all edges $E$ adjacent to $v$,
and the sign of $u(E)$ is chosen with the help of $v$.

\begin{figure}[h]
\includegraphics[height=50mm]{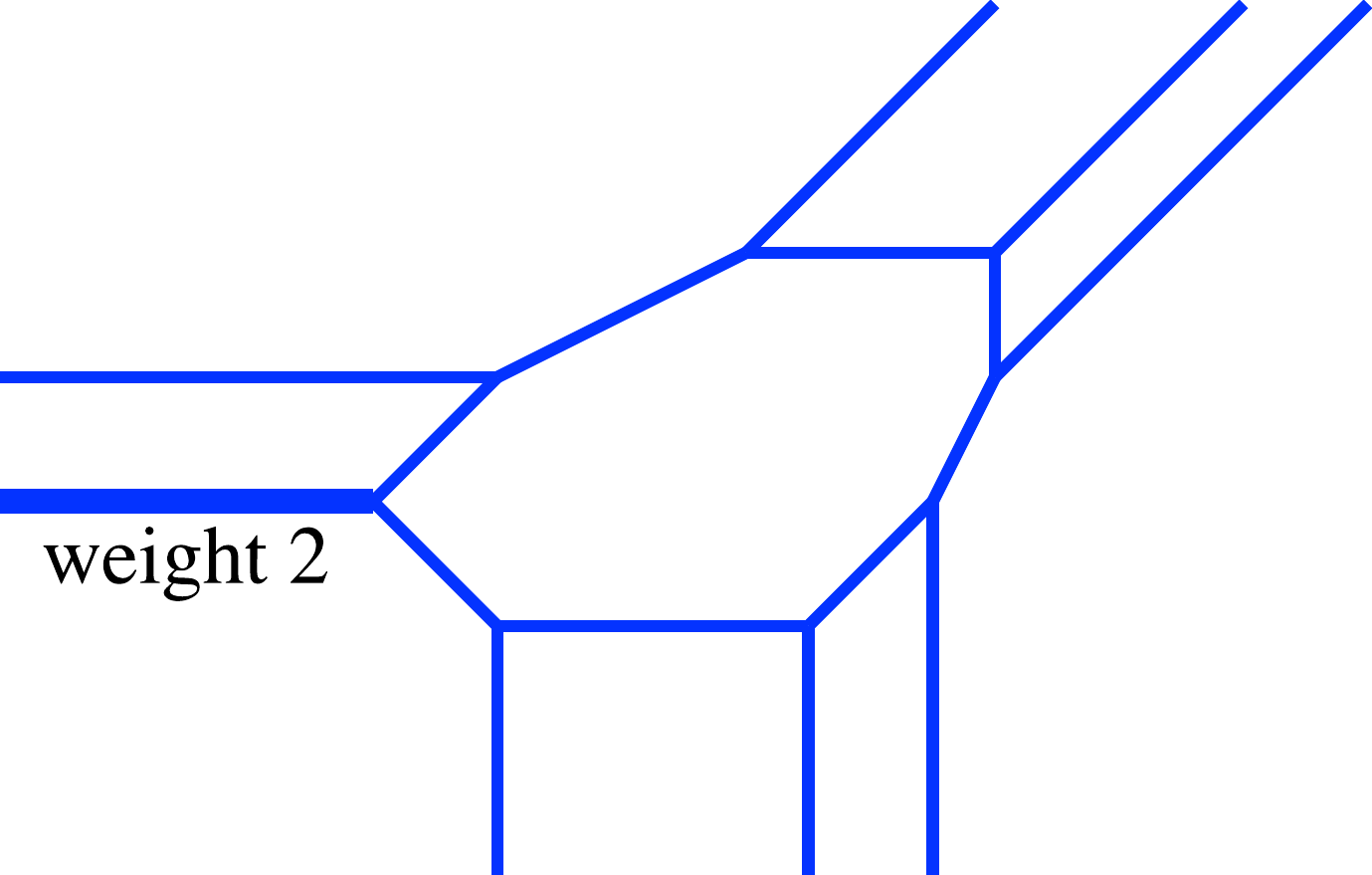}
\caption{\label{cubic} A cubic curve with an edge of weigth 2.}
\end{figure}

Figure \ref{cubic} depicts the corner locus $\Gamma_f$ of a tropical cubic polynomial
$$``f(x,y)=a_{00}+a_{10}x+a_{20}x^2+a_{30}x^3+a_{01}y+a_{11}xy+a_{21}x^2y +
a_{02}y^2+a_{12}xy^2+a_{03}y^3"$$
for some values of $a_{ij}\in\T$.
We leave it as the exercise to the reader to identify the components of $\R^2\setminus\Gamma_f$
with the corresponding monomials. Note that the corner locus
here determines the coefficients $a_{ij}$ up to simultaneous tropical multiplication by a constant.
\footnote{This is due to the fact that every monomial of $f$ in this example
corresponds to some non-empty component of
$\R^2\setminus\Gamma_f$ where it is strictly larger than other monomials.
In general some monomials might be nowhere dominating. Their coefficients are not
determined by $\Gamma_f$.} 

As in classical geometry the same curve can appear inside a plane (or, more
generally, a higher-dimensional space)
in several possible ways. Thus it is useful to define the curve in intrinsic terms,
without referring to the ambient space.
{\em Abstract tropical curves}
are
so-called ``metric graphs".
In the compact case these are
finite connected graphs equipped
with an inner metric such that all edges adjacent to $1$-valent
vertices have infinite length.
More generally, a tropical curve is obtained from such finite graph by removing
some of its 1-valent vertices. 
The complement of all remaining 1-valent vertices
is a metric space. Curves
are considered isomorphic if they are homeomorphic so that the homeomorphism
preserves this metric.

Tropical curves are counterparts
of Riemann surfaces. The role of the genus is played
by the first Betti number ({\em i.e.}, the number of independent
cycles) of the graph.
The role of the punctures is played by the removed 1-valent vertices.
Compact (or projective) tropical curves are finite graphs themselves:
not a single vertex is removed.

Let $C$ be a tropical curve and $x\in C$ be a point which is not a 1-valent vertex.
We may form
a new graph $\tilde C$ from the disjoint union of $C$ and the infinite ray $[0,+\infty]$
(considered as a metric space after removing $+\infty$) by identifying $x$ and $0$.
The result is a compact tropical curve of the same genus and with the same number of
punctures. Furthermore we get a natural contraction map $\tau_x:\tilde{C}\to C$.
The map $\tau_x$ is called {\em tropical modification} at $x$. Tropical modifications generate
an equivalence relation on tropical curves. Any edge connecting a 1-valent vertex and a vertex of valence at least 3
can be contracted.

\begin{figure}[h]
\includegraphics[height=30mm]{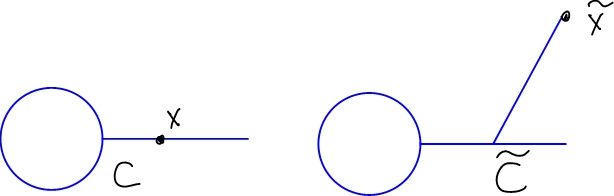}
\caption{\label{modif} Tropical modification}
\end{figure}

\begin{figure}[h]
\includegraphics[height=30mm]{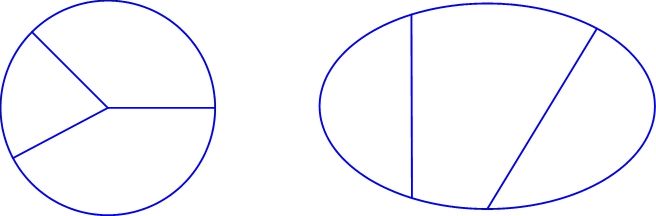}
\caption{\label{curves} Tropical curves of genus 3}
\end{figure}

We arrive to our first example of correspondence between tropical and classical geometric objects.
{\em Compact Riemann surfaces (complex curves) correspond to metric graphs up to tropical modifications (tropical curves).}
A tropical curve of positive genus has a natural {\em minimal model} with respect to tropical modifications.
It is obtained by contracting all edges adjacent to $1$-valent vertices.
Figure \ref{curves} depicts some tropical curves of genus $3$.

It is easy to note that the dimension of the space of tropical curves of genus $g > 1$ is $3g-3$ and thus coincides
with the dimension of the space of complex curves. Most classical theorems on Riemann surfaces have
their tropical counterparts.

We can modify a previous example by marking a number of 1-valent vertices
on a tropical curve.
{\em Riemann surfaces with marked points correspond to metric graphs
with marked points.}
Once a 1-valent vertex is marked it can no longer be contracted
by tropical modifications.
Once at least two points on a rational (genus 0) tropical curve
are marked it also admits a natural
minimal model.

\begin{figure}[h]
\includegraphics[height=30mm]{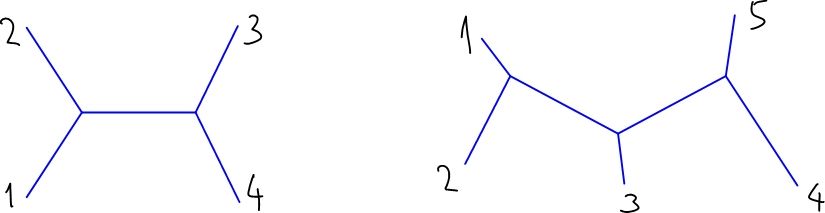}
\caption{\label{mcurves} Rational curves with marked points}
\end{figure}



The only compact tropical higher-dimensional space we consider
in this section is
the tropical projective $n$-space
$$\tp^n = \{(x_0, \ldots x_n) \in (\T^{n + 1} \setminus
\{(-\infty, \ldots, -\infty)\})\}/\sim,$$
where the equivalence relation $\sim$ is defined as follows:
$$(x_0, \ldots x_n) \sim (x'_0, \ldots, x'_n)$$
if and only if
there exists a real number~$\lambda$ such that
$x_i = "\lambda x'_i"$ ({\em i.e.}, $x_i = \lambda + x'_i$)
for any $i = 0$, $\ldots$, $n$.
If $x_0\neq-\infty$ we may take $(x_1-x_0,\dots,x_n-x_0)$ as
affine coordinates, so $\tp^n\setminus\{x_0=-\infty\}=\T^n$ as
in the classical case.
The set defined by $x_j\neq -\infty$, $j=0,\dots,n$,
is $(\T^\times)^n=\R^n\subset\T^n\subset\tp^n$.
This is the finite part of $\tp^n$.

Topologically we may think of $\tp^n$
as an $n$-dimensional simplex. In particular, we get $$\R^n=\tp^n\setminus\dd\tp^n.$$
Tropical structure on each (relatively) open $k$-dimensional face of
$\tp^n$
is a tautological
integer-affine structure on $\R^k$.

This gives another example of the tropical correspondence principle: the complex projective
space $\cp^n$ becomes the $n$-simplex $\tp^n$.

\begin{figure}[h]
\includegraphics[height=30mm]{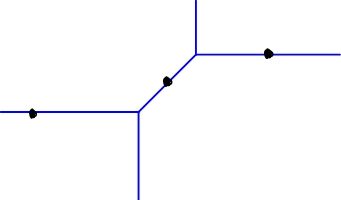}
\caption{\label{tp2} A rational curve in $\tp^2$}
\end{figure}

Up to tropical modifications all compact tropical curves can be embedded in $\tp^n$ by {\em tropical maps},
which are the degenerations in the tropical limit
of holomorphic embeddings in $\cp^n$ of Riemann surfaces.
A tropical map $h:C\to\tp^n$ is a continuous map with the following properties.
\begin{itemize}
\item For every  edge $E\subset C$,
we have $h(E)\subset\R^n$, the map $h|_E$ is
smooth,
and $(d(h|_E))_x(u)\in\Z^n$ at every point $x\in E$ whenever $u$ is a tangent vector of unit length. Note that this condition implies that
$h(E)$ is a straight (possibly unbounded) interval in $\R^n$ whenever $h|_E$ is non-constant.
By continuity, there are only two possible values for $(d(h|_E))_x(u)$ which differ by sign.
Once an endpoint $v$ of the edge $E$ is chosen,
we define $u(E)$ to be equal
to $(d(h|_E))_x(u)$,
where $x$ is a point of $E$, and $u$ is the tangent unit
vector oriented away from $v$. 
The GCD of the absolute values of the components of
$u(E)$ is the {\em weight} $w(E)$ of $E$.
\item For every vertex $v\in C$ we have the balancing condition \eqref{balancing}.
\end{itemize}

Note that we have $h(x)\in\R^n$ for every $x\in C$ unless $x$ is a 1-valent vertex. If $x$ is a 1-valent vertex
then $h(x)\in\R^n$ if and only if the edge adjacent to $x$ has weight $0$, otherwise a 1-valent vertex $x$ is mapped to $\dd\tp^n$.
The set $h(C)\cap\R^n$ is a piecewise-linear graph that can be naturally enhanced with the weights.
The inverse image of an edge $E'\subset h(C)$ may be contained in several
edges $E_1,\dots,E_k$ of $C$. We set
$$w(E')=\sum\limits_{j=1}^k w(E_j).$$

For the case $n=2$ it is a straightforward exercise to see that $h(C)\cap\R^2$ can be presented
as the corner locus $\Gamma_f$ of some tropical polynomial $f$.
For example, after doing further modifications on abstract
tropical elliptic curves depicted on Figure \ref{modif} we can map them to $\tp^2$
so that their image will be presented as the corner locus of a tropical cubic polynomial
as the one from Figure \ref{cubic}. Clearly, the number of 1-valent vertices after modifications
has to agree with the number of ends on the planar picture.

Note that the condition we impose on tropical map implies that the length of the circle in the
metric of abstract tropical curve and the length of the cycle of the planar curve as the one on Figure \ref{cubic}
have to agree. The metric on the planar curve is defined by the condition that $u(E)$ is a unit vector.
In particular, this means that the edge length coincides with the one given by
the Euclidean metric for vertical and horizontal edges of weight $1$.
It is shorter by factor of $\sqrt {2}$ for
diagonal edges of weight $1$. For edges $E$ 
of higher weight we have to additionally divide the length by $w(E)$. 

\begin{figure}[h]
\includegraphics[height=50mm]{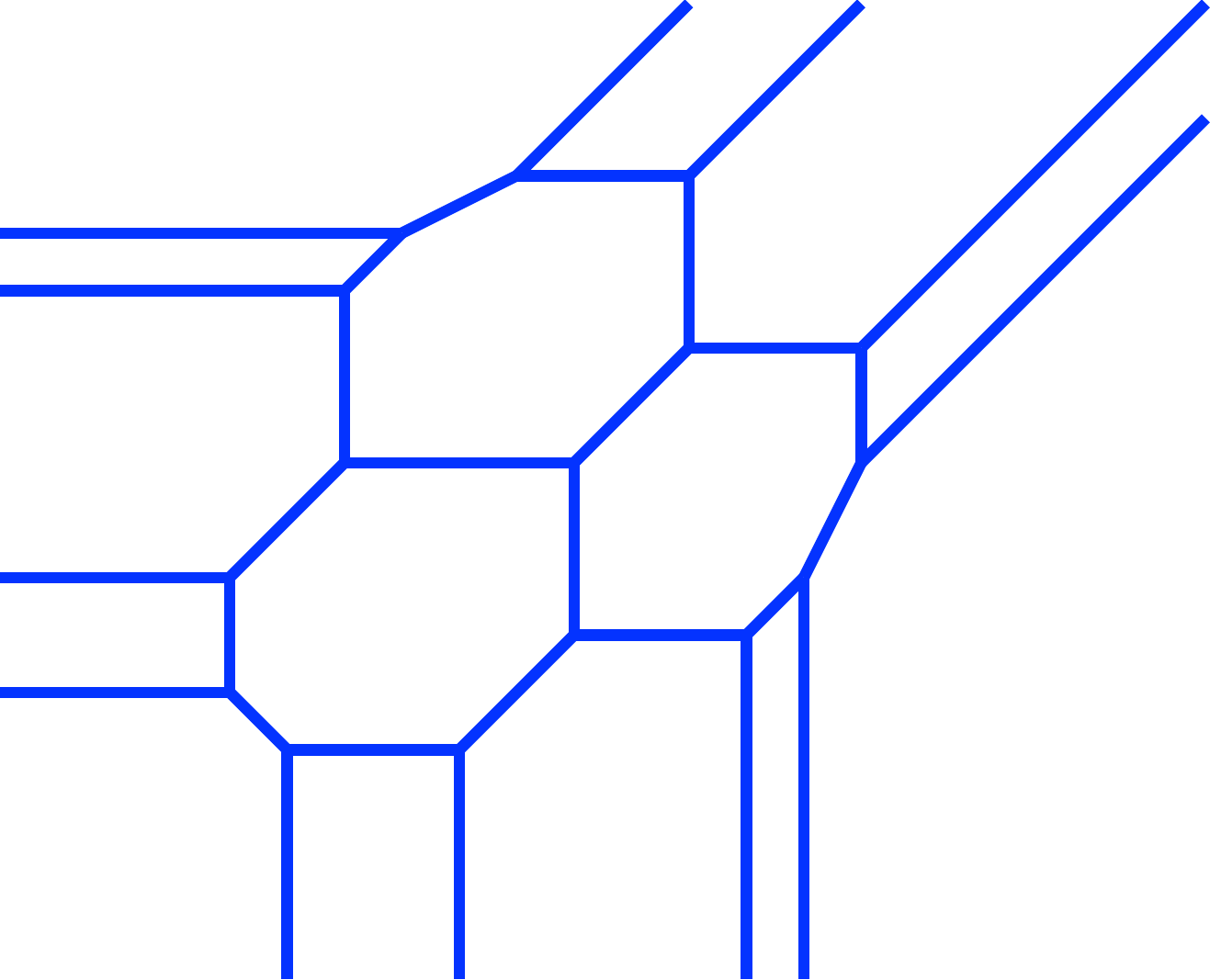}
\caption{\label{quartic} A quartic curve.}
\end{figure}

Figure \ref{quartic} shows a possible image of the left-hand curve from Figure \ref{curves}
under a tropical map to $\tp^2$ after doing 12 tropical modifications.
This image
can be given by 
a quartic tropical polynomial in two variables. Here the length of all three independent circles
have to agree. As in the classical case one can show that any curve of genus 3 can be presented
by a quartic planar curve (up to modifications) unless the curve is {\em hyperelliptic},
{\it i.e.},
admits an isometric involution such that its quotient space is a tree. 

The examples of correspondence
principle
that we considered so far
can be combined to a correspondence between projective complex and projective
tropical curves.
Such a correspondence can be used in applications
to enumerative geometry
(as it was shown in a series of works starting
from Mikhalkin's work~\cite{M}
on
tropical enumerative geometry in $\R^2$).

Tropical approach provides heuristics for many
problems in
classical algebraic geometry (including as it was recently noted by
Kontsevich such a central open problem as the Hodge conjecture).
Each instance of the tropical correspondence is a separate theorem.
Expanding tropical correspondence is an active topic of current research.

\section{Floor diagrams}

The correspondence principle mentioned in the previous section allows one
to reduce certain enumerative problems concerning complex
curves to tropical enumerative problems. How
can we solve the resulting
tropical problems? For example, how
can we enumerate
tropical curves (counted with the multiplicities
dictated by the correspondence) of degree~$d$ and genus~$g$ which pass
through $3d - 1 + g$ points in general position in $\tp^2$?
One of the possible ways of enumeration of tropical curves
is provided by floor diagrams \cite{BM1,FM}.

Choose one of the vertices of the coordinate system in $\tp^2$,
for example, the point $[0: 1: 0]$.
The straight lines which pass through the chosen vertex
and do not pass through any other vertex of the coordinate system
are called {\it vertical}.
Let $T$ be a tropical curve in $\tp^2$.
An edge of~$T$ is called an {\it elevator}
if it is contained in a vertical straight line.
Denote by $El(T)$ the union of elevators and
adjacent vertices of~$T$.
A {\it floor} of~$T$
is a connected component of the closure of the complement
of $El(T)$ in $T$.

Choose now $3d - 1 + g$ points in general position
in $\tp^2$ and ``stretch" the chosen configuration
of points in the vertical direction,
that is, move the points of the configuration
along vertical straight lines in such a way that
the distance between any two points of the configuration
becomes very big (for any two points of the configuration,
one point becomes much ``higher" than the other one).
Denote the resulting configuration by $\omega$.

It is not difficult to check that if a tropical curve of degree~$d$
and genus~$g$
is traced through the points of $\omega$, then
\begin{itemize}
\item the curve contains exactly $d$ floors,
$d - 1 + g$ elevators of finite length,
and $d$ elevators of infinite length
(the latter elevators are adjacent to one-valent vertices
on the coordinate axis $x_1 = -\infty$),
\item each floor and each elevator
of the curve contains exactly one
point of $\omega$.
\end{itemize}

Such a tropical curve can be represented by a connected graph
whose vertices correspond to the floors of the curve and whose edges
correspond to the elevators. This graph
is naturally oriented: each elevator of the tropical
curve can be directed toward the point $[0: 1: 0]$,
{\it i.e.}, vertically up.

A {\it floor diagram} of degree~$d$ and genus~$g$ is
a connected oriented weighted (each edge has a positive
integer weight) graph $D$ such that
\begin{itemize}
\item the graph~$D$ is acyclic as an oriented graph,
\item the first Betti number $b_1(D)$ of~$D$ is equal to~$g$,
\item the graph~$D$ has exactly $d$ {\it sources}, that is, one-valent
vertices whose only adjacent edge is outgoing,
\item any edge adjacent to a source is of weight $1$,
\item for any vertex~$v$ of~$D$ such that $v$ is not a source,
the difference between the total weight of ingoing edges of~$v$
and the total weight of outgoing edges of~$v$ is equal to~$1$.
\end{itemize}

Each floor diagram of degree~$d$ and genus~$g$ has
$2d$ vertices ($d$ of them are sources and $d$ others are not sources)
and $2d - 1 + g$ edges.
Denote by $M(D)$ the union of the set of edges of~$D$ and
the set of vertices of~$D$ which are not sources.
The set $M(D)$ is partially ordered.
We say that a map $m$ between two partially ordered sets
is {\it increasing} if $m(i) > m(j)$ implies $i > j$.
A {\it marking} of a floor diagram $D$ of degree~$d$ and genus~$g$ is an increasing bijection $m: \{1, 2, \ldots, 3d - 1 + g\} \to M(D)$.
A floor diagram equipped with a marking is called
a {\it marked floor diagram}.

Assume that the points of the configuration~$\omega$ considered above
are numbered by the elements of $\{1, 2, \ldots, 3d - 1 + g\}$
in the increasing order of heights of the points.
Then, any tropical curve of degree~$d$ and genus~$g$
which passes through the points of~$\omega$ gives
rise to a marked floor diagram of degree~$d$ and genus~$g$.
Reciprocally, any marked floor diagram
of degree~$d$ and genus~$g$ gives rise to a tropical curve
of degree~$d$ and genus~$g$ which passes through the points of~$\omega$,
so we get a 1-1 correspondence between marked floor diagrams and tropical
curves passing through $\omega$.

\begin{figure}[h]
\includegraphics[height=50mm]{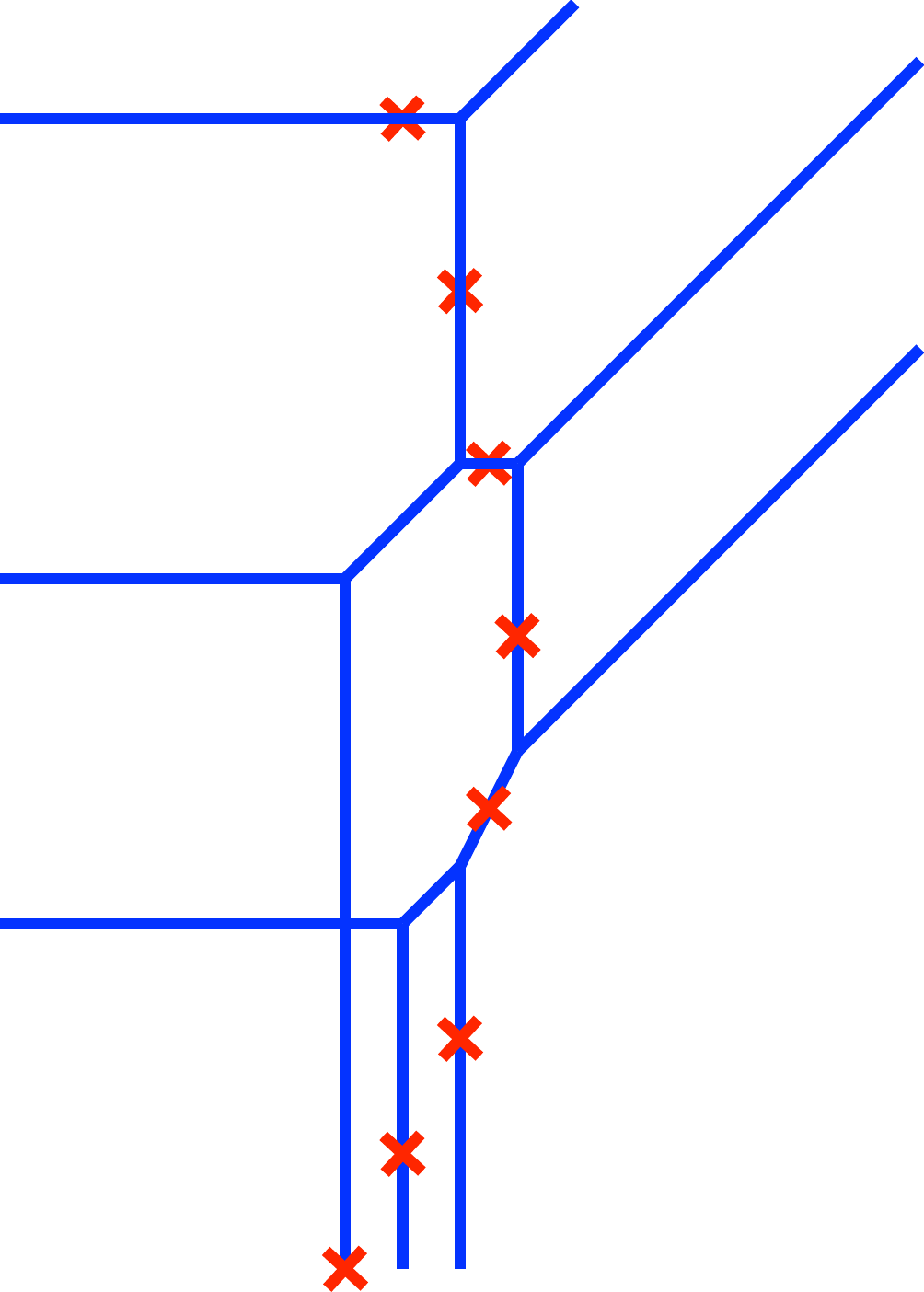}
\caption{\label{reconstruction} A cubic curve passing 
through $8$ vertically stretched points.}
\end{figure}

Figure \ref{reconstruction} shows the tropical curve corresponding to the first marked floor
diagram from Figure \ref{3d3-1} for a choice of a generic vertically stretched configuration $\omega$ of 8 points. 
We leave it as an exercise to the reader to reconstruct the tropical curve corresponding to other
marked floor diagrams for the same choice of the configuration $\omega$. 

Thus, to enumerate the tropical curves (counted with the multiplicities
dictated by the correspondence) of degree~$d$ and genus~$g$ which pass
the points of~$\omega$, it is enough to enumerate the marked floor
diagrams (counted with appropriate multiplicities) of degree~$d$
and genus~$g$. It turns out that, for any marked floor diagram,
the appropriate multiplicity to consider is the product
of squares of weights of the edges. By \cite{M} the sum of multiplicities
of all marked diagrams of degree $d$ and genus $g$ with these multiplicities is
equal to the number of all curves of degree $d$ and genus $g$ passing through
a configuration of $3d-1+g$ generic points in $\cp^2$.

\begin{exa}
To compute the number of rational cubic curves  passing
through 8 generic points in $\cp^2$ we need to enumerate
marked floor diagrams of genus 0 with 3 sources.
Before marking there are only 3 such diagrams, see Figure \ref{3d3}.

\begin{figure}[h]
\includegraphics[height=45mm]{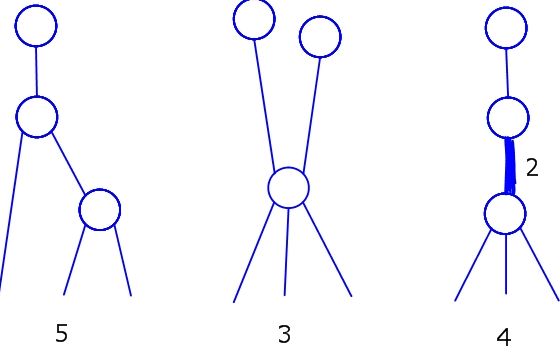}
\caption{\label{3d3} Floor diagrams enumerating rational cubic curves 
in the plane}
\end{figure} 

Here the vertices of the diagrams other than sources are shown with small circles.
All sources are placed in the bottom of the diagrams. Each edge is oriented upwards.

The first diagram supports five different markings, see Figure \ref{3d3-1}, the second one support three different
markings, see Figure \ref{3d3-2}.
\begin{figure}[h]
\includegraphics[height=30mm]{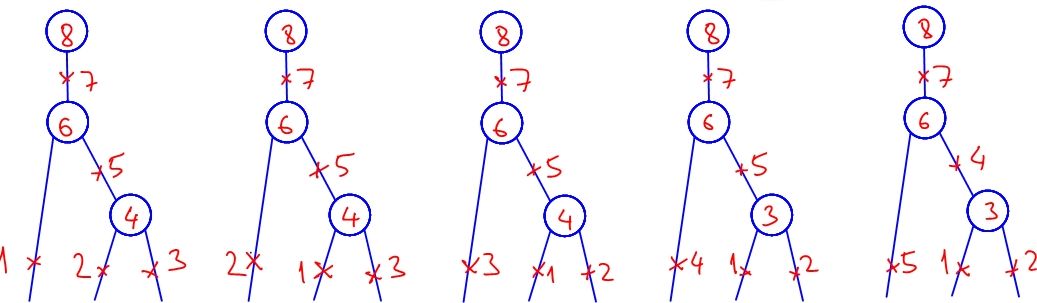}
\caption{\label{3d3-1} Markings for the first diagram in Figure \ref{3d3}}
\end{figure}
\begin{figure}[h]
\includegraphics[height=30mm]{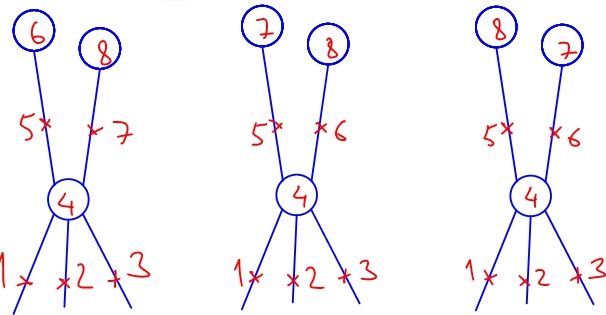}
\caption{\label{3d3-2} Markings for the second diagram in Figure \ref{3d3}}
\end{figure}
The last one supports only one marking, but comes
with multiplicity 4 as it contains a weight 2 edge.
Adding $5+3+4$ we get 12 rational cubic curves passing through 8 generic points in $\cp^2$.
\end{exa}

\begin{exa}
To consider a more complicated example we consider the problem of enumeration of
degree 4 curves of genus 1 in $\cp^2$. We get 11 diagrams before we take marking in consideration.
Figure \ref{4d} indicates the number of markings taken with multiplicities.
As the result we get $26+16+15+24+9+9+21+28+21+32+24=225$ elliptic quartic curves through
12 generic points in $\cp^2$.
\begin{figure}[h]
\includegraphics[height=70mm]{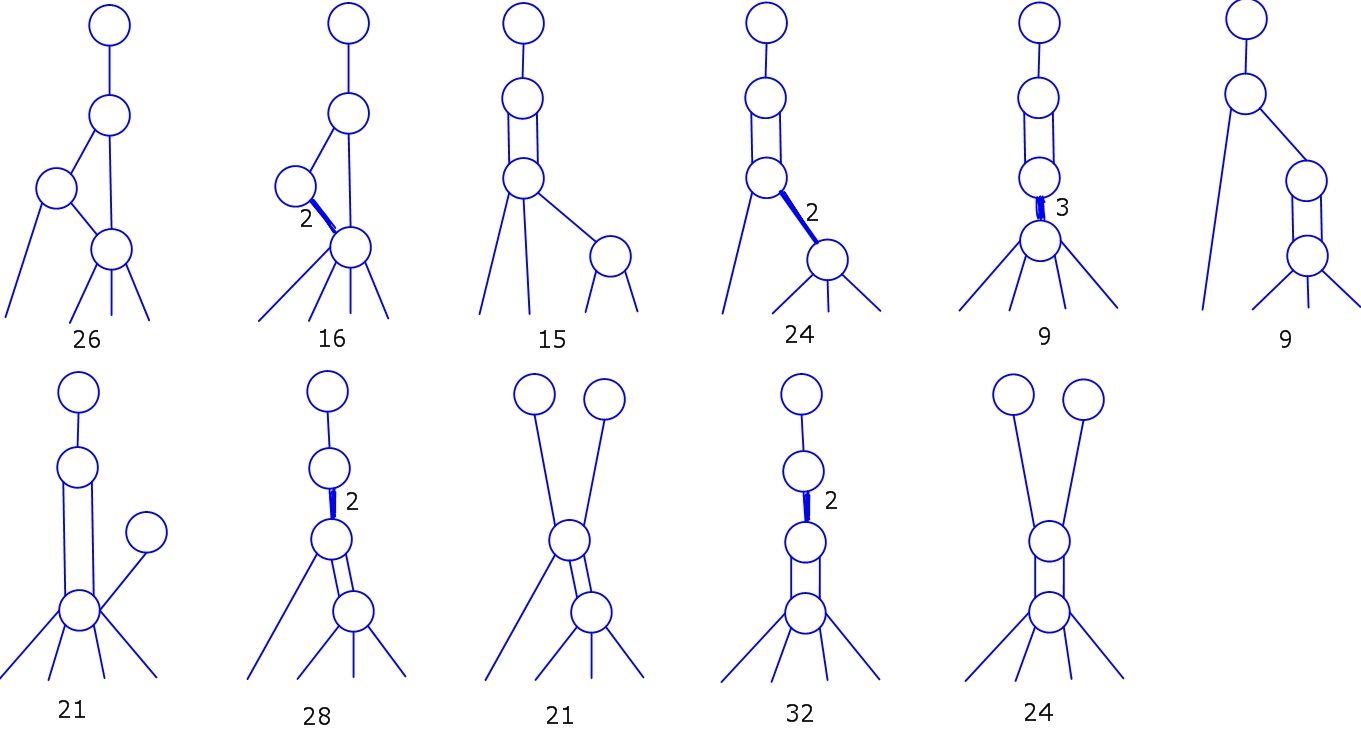}
\caption{\label{4d} Floor diagrams and numbers of their markings (with multiplicities) for the degree 4, genus 1 case}
\end{figure}
\end{exa}

\begin{rmk}
There is also a correspondence between floor diagrams and real algebraic curves of
degree $d$
and genus $g$
which pass
through appropriately chosen $3d-1+g$ points
in general position in $\rp^2$. 
We can introduce the {\em real multiplicity} of a floor diagram to be zero if the diagram has an edge of
even weight and 1 otherwise. Denote by $N_{\R}(g,d)$ the sum
of
real multiplicities over all floor diagrams
of degree $d$ and genus $g$.
Computing real multiplicities in the examples above gives us
$N_\R(0, 3) = 8$
and
$N_\R(1, 4) = 93$.

It turns out that there always
exists
a configuration of $3d - 1 + g$ generic points in $\rp^2$
so that
there are
at least $N_{\R}(g,d)$ real curves
of degree $d$ and genus $g$
passing through them.
These real curves are nodal, and a real node
of a real curve can either be {\em hyperbolic}
(an intersection of two real branches of the curve) or {\em elliptic} (an intersection
of two conjugate imaginary branches of the curve). Denote the number of elliptic nodes by $e$.
If we enhance each real curve with the sign $(-1)^e$ as suggested by Welschinger \cite{W}, then
the corresponding number of all real curves
of degree $d$ and genus $g$
through
our configuration will be equal to $N_{\R}(g,d)$, see \cite{M}.

In \cite{W} it was shown that
the number of real curves, counted with signs $(-1)^e$,
of degree $d$ and genus $g$
which pass through $3d - 1 + g$ points in $\rp^2$
does not depend on the choice of the configuration of
points
as long as this configuration is generic
and $g = 0$.
An interesting phenomenon occurs for $g>0$:
this number is {\em not} invariant in the context of classical real algebraic geometry, but it {\em is} invariant
in the context of tropical geometry (see~\cite{IKS2}).
This area is currently a subject of active research,
see relevant discussions in \cite{IKS1}, \cite{IKS2} and \cite{M-mfo}.
\end{rmk}

\end{document}